\documentclass[twocolumn]{amsart}

\usepackage{amssymb,amsmath,amscd,graphicx}
\usepackage{caption,array,calc}
\usepackage[export]{adjustbox}

\usepackage{subfig}
\newlength\celldim \newlength\fontheight \newlength\extraheight
\newcounter{sqcolumns}
\newcolumntype{S}{ @{}
>{\centering \rule[-0.5\extraheight]{0pt}{\fontheight + \extraheight}}
p{\celldim} @{} }
\newcolumntype{Z}{ @{} >{\centering} p{\celldim} @{} }
\newenvironment{squarecells}[1]
{\setlength\celldim{1.1cm}%
\settoheight\fontheight{A}%
\setlength\extraheight{\celldim - \fontheight}%
\setcounter{sqcolumns}{#1 - 1}%
\begin{tabular}{|S|*{\value{sqcolumns}}{Z|}}\hline}
{\end{tabular}}
\newcommand\nl{\tabularnewline\hline}

\usepackage{tikz}
\tikzset{every picture/.style={line width=0.75pt}} %set default line width to 0.75pt
\usetikzlibrary{shapes.geometric}
\usetikzlibrary{quotes,angles,positioning}

\DeclareMathAlphabet{\mymathbb}{U}{bbold}{m}{n}

\newcommand{\ts}{\hspace{0.5pt}}

\newtheorem{theorem}{Theorem}

\theoremstyle{definition}

\newcommand{\ee}{\ts\mathrm{e}}
\newcommand{\dd}{\,\mathrm{d}\ts}

\newcommand{\dens}{\mathrm{dens}}

\newcommand{\ii}{\ts\mathrm{i}}

\renewcommand{\epsilon}{\varepsilon}

\begin{document}	
	\twocolumn[
	\begin{LARGE}
		\centerline{Substitutions and their generalisations}\vspace{3ex}
	\end{LARGE}
	\centerline{\large Neil Ma\~nibo$^{1}$}\vspace{2ex}
	\begin{footnotesize}
		\centerline{
			${}^{1}$\textit{Bielefeld University, Bielefeld, Germany \vspace{1mm}}}
	\end{footnotesize}\vspace{4ex}
	\begin{small}
		\hrule\vspace{2ex}
		\begin{minipage}{\textwidth}
			\textbf{Abstract}\vspace{2ex}\\ 
Tilings and point sets arising from substitutions are classical mathematical models of quasicrystals. Their hierarchical structure allows one to obtain concrete answers regarding spectral questions tied to the underlying  measures and potentials. In this review, we present some generalisations of substitutions, with a focus on substitutions on compact alphabets, and with an outlook towards their spectral theory. Guided by two main examples, we will illustrate what changes when one moves from finite to compact (infinite) alphabets, and discuss under which assumptions do we recover the usual geometric and statistical properties which make them viable models of materials with almost periodic order. 
We also present a planar example (which is a two-dimensional generalisation of the Thue--Morse substitution), whose diffraction is purely singular continuous. 
			\vspace{2ex}\\
			\textit{Keywords:}\/  substitutions, tilings, diffraction, Schr\"odinger operators, compact alphabets
		\end{minipage}\vspace{2ex}
		\hrule
	\end{small}\vspace{6ex}
	]

	\section{Introduction }\label{sec:intro}

Substitutions are well studied objects in the theory of aperiodic order. In symbolic dynamics, they give rise to prototypical models of infinite shift spaces with nice properties. In tiling theory, they generate tilings of Euclidean space, whose aperiodicity is easily confirmed. This is due to the presence of hierarchical structures which are incommensurate with having a non-trivial lattice of translational periods. Both the well-known Fibonacci tiling and Penrose tiling can be generated by a substitution rule. 

Within dynamical systems theory, there is a rich literature on substitutions going back from the 60s to the late 70s; see \cite{Martin, Michel,Dekking}. After the discovery of quasicrystals in 1982, more detailed investigations of their spectral theory were carried out; see \cite{BomTaylor} and \cite{GodLuck}, for instance.

We refer the reader to the monographs by Baake and Grimm \cite{BGr} and N. Pytheas Fogg \cite{P-Fogg}
for detailed expositions on substitutions. We recommend the monographs  by Allouche and Shallit \cite{AS} for connections to theoretical computer science, Kellendonk, Lenz, and Savinien, and Akiyama and Arnoux \cite{KLS,AA} for several articles on different directions (Schr\"odinger operators, algebraic topology, dynamical systems and number theory) and Queff\'{e}lec \cite{Queffelec} for a detailed treatment of the dynamical spectrum.

\begin{figure}[!h]
\subfloat[The GLB inflation rule]{
\includegraphics[scale=0.7]{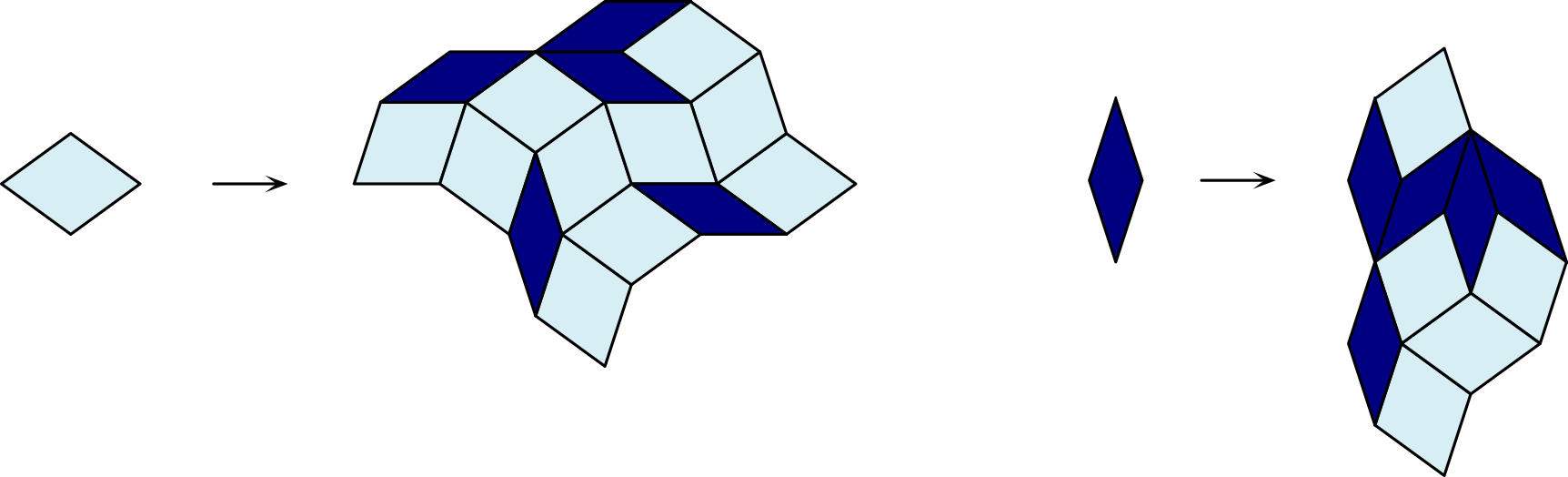}}

\subfloat[A finite patch of the GLB tiling. Taken from the Tilings Encyclopedia \cite{FGH}.]{
\includegraphics[scale=0.85]{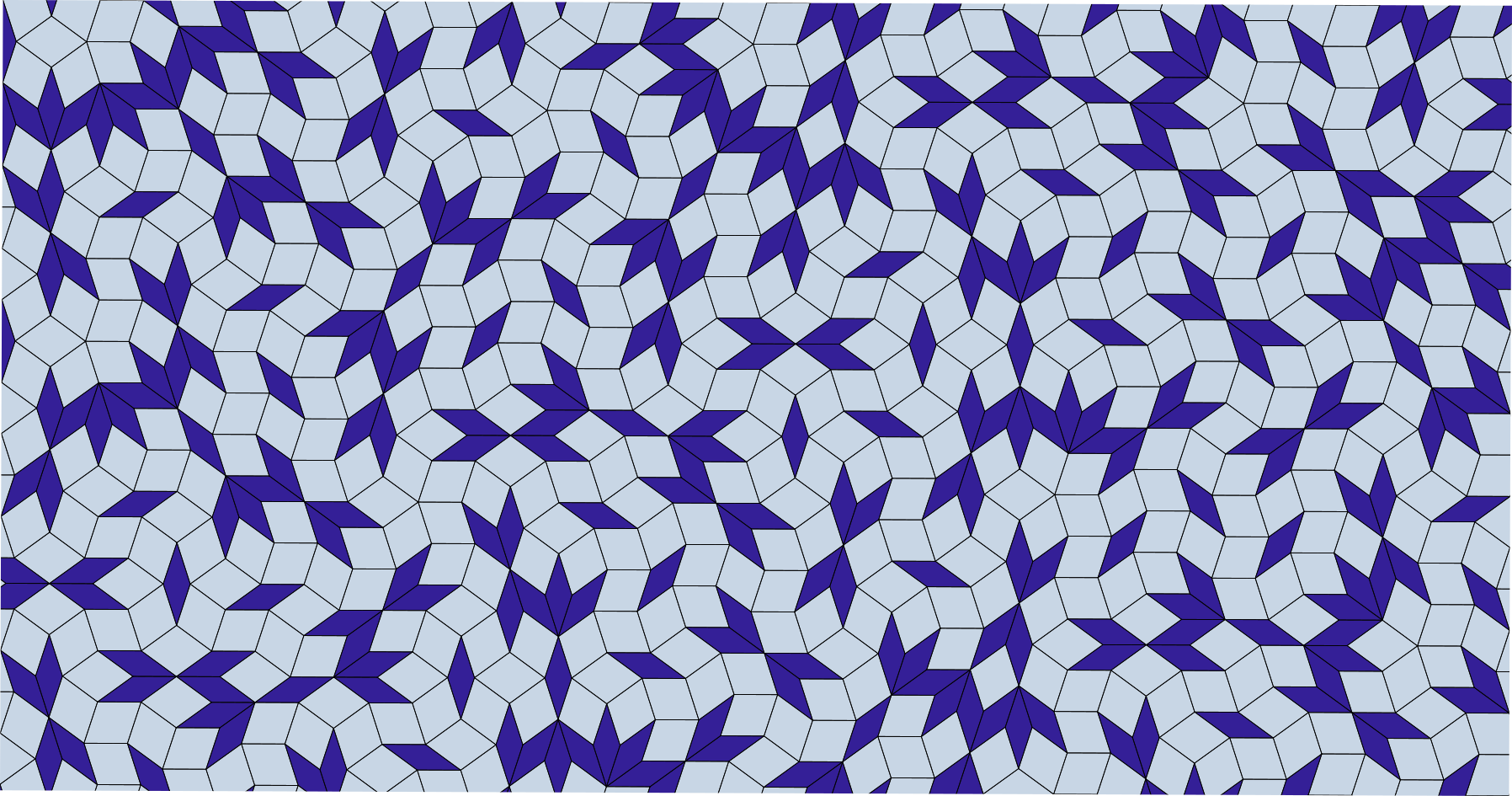}
}
\caption{The GLB tiling is the first two-dimensional tiling which is proved to have essentially purely singular continuous diffraction \cite{BGM}.}\label{fig:GLB}
\end{figure}

Let say a few more words on recent developments regarding the spectral theory of substititution dynamical systems. The \emph{pure point} component is relatively well understood, where new results linked to almost periodicity provide a complete characterisation of pure point diffraction at the level of the atomic configuration; see \cite{MStr-AO} for a survey on almost periodic measures and \cite{Strungaru, LSS-2,LSS-3}. 
In the Pisot case, the situation is even better with the availability of the cut-and-project formalism, which gives an alternative viewpoint with which one can investigate the underlying point sets; see \cite{BGr, Moody-1, Moody-2, Strungaru-2} and \cite{Sing} also for the treatment of the non-unit case. 

The study of the \emph{continuous} component of the diffraction (which corresponds to \emph{diffuse scattering} in experiments) is more challenging, but has seen remarkable progress in the last 10 years alone. A cornerstone of the theory relies on renormalisation techniques, which one can transfer from the level of point sets/tilings to the underlying measures. This allows for a separate treatment of the three (mutually singular) spectral components (pure point, absolutely continuous, and singular continuous). 

To rule out the presence of an absolutely continuous component, a powerful tool is the \emph{Fourier cocycle/spectral cocycle} and the corresponding Lyapunov exponent; see \cite{BFGR, BGM,BufSol, BufSol2, SolTrevino}. Other singularity results for substitutions have also been developed; see \cite{Bartlett,BerSol}.

In particular, one can use this to prove that the Godr\`{e}che--Lan\c{c}on--Billard (GLB) tiling \cite{GodLancon} has a purely singular continuous diffraction apart from the trivial Bragg peak at zero; see \cite{BGM,Manibo} for details. The inflation rule and a patch of the tiling are shown in Figure~\ref{fig:GLB}.

For some systems with singular continuous diffraction, a finer description of the spectrum can be obtained via multifractal analysis; see \cite{BGKS, GLS, FSS} and references therein. A classical example is the Thue--Morse substitution, which we will revisit later.

From a bi-infinite point $w$ generated by a substitution, one can also define the corresponding one-dimensional (discrete) Schr\"odinger operator. These are well studied objects, and aperiodic sequences give rise to aperiodic potentials which inherit the hierarchical structure of the substitution; see the recent monograph \cite{DF} by Damanik and Fillman for the general theory in one dimension.  A second volume \cite{DF-2} on specific classes (including quasiperiodic potentials) is currently in preparation.  For a survey article on Schr\"odinger operators in the context of mathematical quasicrystals, see \cite{DEG}.

In particular, it is well known that the spectrum of such operators generated by aperiodic primitive  substitutions on finite alphabets is always a Cantor set of Lebesgue measure zero. Moreover, having purely singular continuous spectral measures is also common (as opposed to having eigenvalues with exponentially-decaying eigenfunctions, which is typical in random systems, or having absolutely continuous spectral measures which is the case for periodic potentials); see \cite{DGY} for an account on the Fibonacci Hamiltonian.

From these two spectral viewpoints, it is clear that the presence of a well-defined hierarchical structure is reflected in the spectrum. Moreover, it also makes statistical computations (that is, relative frequencies of patches, orbit- and space-averages of certain functions, pair correlations, etc.) tractable and explicit (using tools from ergodic theory). 

One disadvantage of using substitutions on finite alphabets or finitely many prototiles is that they are too rigid. A common consequence (although not guaranteed) of working with substitutions with finitely many starting blocks (letters or prototiles) is that the infinite structure generated has \emph{finite local complexity} (FLC). This means, given a fixed volume $R>0$, there are only finitely many patches of this given volume.  From a computational and theoretical standpoint, this is satisfactory, as one normally deals with structures which admit finitely many local configurations (such as clusters). On the other hand, this usually creates a disparity from what one sees experimentally, as atoms in materials are rarely stationary. 

\begin{figure}[!h]
\begin{center}
\includegraphics[scale=0.8]{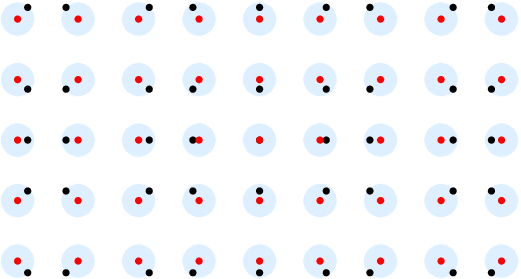}
\end{center}
\caption{A patch from the modulated lattice $\varLambda_f=\left\{(m,n)+f(m,n)\colon (m,n)\in \mathbb{Z}^2\right\}$ (black) with  $f(m,n)=(\frac{1}{5}\sin(2\pi m\sqrt{2}),\frac{1}{5}\sin(2\pi n\sqrt{7}))$. The point set $\varLambda_f$ is Delone but has ILC. Nevertheless, it has pure point diffraction. The lattice $\mathbb{Z}^2$ (red) and the range of deviation from a lattice point (blue disk) are added for visualisation.  
}\label{fig:mod-lat}
\end{figure}

It is thus not surprising that FLC is \emph{not} necessary to have pure point diffraction. As an example, \emph{modulated lattices} normally have infinite local complexity (ILC), but if the modulation is well behaved (for instance, a Bohr almost-periodic function), one still gets pure point diffraction; see Figure~\ref{fig:mod-lat} and \cite{LLRSS} for a recent account on modulated quasicrystals, which includes deformed (weighted) model sets.

Recent developments in the theory of substitution dynamical systems try to circumvent this level of rigidity by considering several generalisations of substitutions while still maintaining the existence of compatible hierarchies (or the semblance thereof). Among these generalisations are substitutions on infinite (yet compact) alphabets; see \cite{DOP,MRW, MRW-2, FGM, FGM-2}. We shall see that, despite dealing with infinitely many symbols (which will be our building blocks), certain natural restrictions allow us to generate geometric model that display properties not possible in the finite setting. At the same time, they possess generalisations of properties in the finite alphabet setting which allow for explicit computations of certain physically-relevant quantities. 
As they give rise to new phenomena, they are certainly interesting from a theoretical standpoint. Likewise, because they give rise to concrete structures with well-defined geometry and statistics, their potential use for actual models of materials (both solid and soft matter) is also promising.

This overview on substitutions on compact alphabets is outlined as follows. In Section 2, we revisit the setting of substitutions on finite alphabets, where we will look at motivating examples and their known dynamical and spectral properties. In Section 3, we  provide a brief survey of generalisations of substitutions which exist in the literature, and point to standard references for the interested reader. In Section 4, we zoom into the world of substitutions on compact alphabets, where we present the generalisation of objects and results associated with substitutions in the finite alphabet setting. 
The consequences of the results discussed and motivated in Section 4 will be more transparent in Section 5, where we  look at several examples which demonstrate the strength and applicability of these results. We  also consider a two-dimensional example and discuss some of its spectral properties. We end with some future directions in Section 6.

\section{Substitutions on finite alphabets revisited}\label{sec:1d-finite}

We start with an example of a substitution on a two-letter alphabet $\mathcal{A}=\left\{a,b\right\}$ given by the rule 
\[
\varrho\colon \begin{cases}
a\mapsto abbb\\
b\mapsto a.
\end{cases}
\]
One can place the (legal) seed $aa$ at the origin (which we denote by $a|a$), and iterate the square $\varrho^2$ of this substitution to obtain a sequence of nested words of growing length. More explicitly, we have 
\[
a|a\overset{\varrho^2}{\longmapsto} abbbaaa|abbbaaa \overset{\varrho^2}{\longmapsto} \cdots \overset{\varrho^2}{\longmapsto} w=\varrho^{2}(w).
\]
Here, $w$ is called a \emph{bi-infinite fixed point} of $\varrho^2$. We consider the (left) shift action $\sigma$ defined on the space of bi-infinite words over $\mathcal{A}$ via $\sigma(x)_n=x_{n+1}$. We then build the orbit-closure of $w$ under $\sigma$, which is given by 
\[
X_{\varrho}=\overline{\left\{\sigma^n(w)\colon n\in \mathbb{Z}\right\}}
\]
where the closure is taken with respect to the \emph{local topology}. In this topology, two bi-infinite words $x$ and $y$ are ``close" if they agree in a large region around the origin. As an example, 
if we consider the words
\begin{align*}
x^{ }_{0}=w&=\cdots abbbaa\textcolor{blue}{a}|\textcolor{blue}{a}bbbaaa\cdots,\\
x^{ }_{1}=\sigma^{1}(w)&=\cdots bbbaaa\textcolor{blue}{a}|\textcolor{red}{b}bbaaaa\cdots,\\
x^{ }_{-1}=\sigma^{-1}w&=\cdots aabbba\textcolor{blue}{a}|\textcolor{blue}{a}abbbaa\cdots,
\end{align*}
one can say that $x^{ }_0$ is \emph{closer} to $x^{ }_{-1}$
than to $x^{ }_{1}$.

The space $X_{\varrho}$ is compact by construction, and consists of all shifted versions of $w$, plus all other bi-infinite words over $\mathcal{A}$ which can be expressed as limits of these shifted words. We call $X_{\varrho}$ the \emph{subshift} associated to $\varrho$. Together with $\sigma$, $(X_{\varrho},\sigma)$ becomes a \emph{topological dynamical system}, where $\sigma$ induces a continuous $\mathbb{Z}$-action on $X_{\varrho}$. We refer the reader to \cite{BGMaz} for a similar exposition focused on the (ramifications of the) Fibonacci substitution. 

%
%This symbolic interlude might seem technical and unnecessary at first. However, as hinted at in the introduction, this dynamical systems approach has the advantage in that statistical calculations on $w$ can be interpreted as calculations on the subshift $X_{\varrho}$.

So far, we have only harvested symbolic sequences from $\varrho$. In what follows, we derive geometric objects (e.g., tilings and point sets) from it. The crucial ingredient to this construction is the \emph{substitution matrix} $M_{\varrho}$ of $\varrho$, which is an $|\mathcal{A}|\times|\mathcal{A}|$-matrix that encodes some combinatorial information of $\varrho$. This is defined entry-wise via 
\[
\left(M_{\varrho}\right)^{ }_{ij}= \text{number of letters of type }i \text{ in } \varrho(j).
\]
For the substitution above, one has 
\[
M_{\varrho}=\begin{pmatrix}
1 & 1\\
3 & 0
\end{pmatrix}.
\]
To see this clearly, note that the first column reads $(1,3)^{T}$ because $\varrho(a)=abbb$ contains one $a$ and three occurrences of $b$. This matrix is \emph{primitive}, that is, it has a power $\left(M_{\varrho}\right)^{n}$ which only has strictly positive entries. If this holds, we also call the  
 substitution $\varrho$  primitive.  Equivalently, for the same power $n$ and for any $a\in\mathcal{A}$, $\varrho^{n}(a)$ contains all letters from $\mathcal{A}$.

Primitivity has an immediate consequence which is the \emph{minimality} of $(X_{\varrho},\sigma)$. 
This is equivalent to every element $x\in X_{\varrho}$ being \emph{repetitive}, meaning 
any finite subword $v$ of $x$ occurs with bounded gaps. 

As a consequence of Perron--Frobenius (PF) theory, primitivity also implies that $M_{\varrho}$ admits a unique eigenvalue of maximum modulus, which we call $\lambda_{\text{PF}}$. Furthermore, this eigenvalue has a one-dimensional left (resp. right) eigenspace, generated by a strictly positive eigenvector  $\boldsymbol{L}$ (resp. $\boldsymbol{R}$).  

For our example, we have 
\[
\lambda=\frac{1+\sqrt{13}}{2}\approx 2.302776.
\]
The second (largest) eigenvalue of $M_{\varrho}$ is $\lambda_2=\frac{1-\sqrt{13}}{2}$, which is an algebraic conjugate of $\lambda$. Since $|\lambda_2|>1$, $\lambda$ is a non-PV number (where PV stands for Pisot--Vijayaraghavan). We then call $\varrho$ a \emph{non-PV substitution}. This example is studied in detail in \cite{BFGR}.

The corresponding left and right PF eigenvectors are given by 
\[
\boldsymbol{L}=(\lambda,1)\quad \text{ and }\quad \boldsymbol{R}=\left(\frac{\lambda-1}{3},\frac{4-\lambda}{3}\right)^{T}.
\] 
The right PF eigenvector $\boldsymbol{R}$ (which is statistically normalised) gives the relative letter frequencies of $a$ and $b$, which are uniform across all elements $x\in X_{\varrho}$. In this case, since the entries of $\boldsymbol{R}$ are both irrational, it follows immediately that every $x\in X_{\varrho}$ is non-periodic (and hence both the substitution and the subshift are called \emph{aperiodic}).
Primitivity implies something stronger in terms of frequencies, namely, that the $(X_{\varrho},\sigma)$ is \emph{uniquely ergodic}. In the language of symbolic dynamics, this is equivalent to the uniform existence of relative frequencies of finite words.

\begin{figure}[!h]
\begin{center}
\includegraphics[scale=0.5]{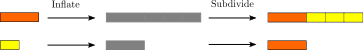}
\caption{The inflation rule associated to the substitution $\varrho\colon a\mapsto abbb, b\mapsto a$. }\label{fig:inflation-rule}
\end{center}
\end{figure}

On the other hand, the left PF eigenvector $\boldsymbol{L}$ allows one to assign a tile $\mathfrak{t}_a$ (which is an interval) of length $\boldsymbol{L}_{a}$ to the letter $a$, for every letter in the alphabet.
This yields a geometric \emph{inflation rule} from $\varrho$ (which by abuse of notation we will also refer to as $\varrho$). 
The image of $\mathfrak{t}_a$ under $\varrho$ can be obtained by first inflating it by a factor of $\lambda$, and then subdividing this blown-up interval into a concatenation of tiles according to the substitution rule; see Figure~\ref{fig:inflation-rule}.

\begin{figure}[!h]
\begin{center}
\includegraphics[scale=0.4]{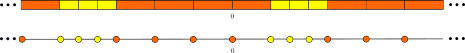}
\caption{A tiling generated by the corresponding inflation rule $\varrho$, together with the derived point set. }\label{fig:inflation-tiling}
\end{center}
\end{figure}

From this inflation rule, one can generate a tiling $\mathcal{T}$ of the line by the tiles $\mathfrak{t}_a$ (orange) and $\mathfrak{t}_b$ (yellow) by placing two copies of $\mathfrak{t}_a$ at the origin and applying the $\varrho$ iteratively. From $\mathcal{T}$, one can obtain a (coloured) point set $\varLambda$ consisting of two ``atom" types by identifying every tile with its left endpoint; see Figure~\ref{fig:inflation-tiling}. This point set is Delone (uniformly discrete and relatively dense) and is aperiodic (via the same frequency argument for the symbolic fixed point).

From $\mathcal{T}$, one can analogously obtain a dynamical system by collecting all $\mathbb{R}$-translates of $\mathcal{T}$ and taking the closure under the \emph{local rubber topology} (which in the case of tilings with FLC, coincides with the local topology); see \cite{BGr}.  One calls
\[
\Omega:=\Omega_{\mathcal{T}}=\overline{\left\{t+\mathcal{T}\colon t\in \mathbb{R}\right\}}
\]
the geometric hull associated to $\mathcal{T}$. As in the symbolic case, $(\Omega,\mathbb{R})$ (equipped with the action of $\mathbb{R}$ by translation) is a topological dynamical system. 
 
We are now poised to discuss some spectral properties of this model. We begin with the diffraction measure. From the point set $\varLambda$ above, one can define a \emph{weighted Dirac comb} 
\[
\mu=\sum_{x\in \varLambda}w(x)\delta_x,
\]
where $\delta_x$ is the Dirac measure at $x$ and $w(x)\in \mathbb{C}$ is the weight, which depends on the atom type at position $x$. From $\mu$, one obtains the \emph{autocorrelation measure} $\gamma^{ }_{\mu}$ given by $\gamma^{ }_{\mu}=\mu\circledast\widetilde{\mu}$, where $\circledast$ is the (volume-averaged) \emph{Eberlein convolution}. This measure is positive definite, and hence Fourier transformable. Its Fourier transform $\widehat{\gamma^{ }_{\mu}}$ is the \emph{diffraction}. For detailed accounts on diffraction of translation-bounded measures, we refer the reader to \cite{BGr, LSS-2, LSS-3}. 
We note that the diffraction does not depend on the choice of $\mathcal{T}\in {\Omega}$ (for a fixed choice of weight function $w$).

For point sets arising from a self-similar tiling with inflation factor $\lambda$,
Solomyak's criterion \cite{Solomyak}, together with Dworkin's argument \cite{Dworkin}, states that the corresponding diffraction measure can only have a non-trivial Bragg peak if $\lambda$ is a PV number. Returning to our example, 
since $\lambda=\frac{1+\sqrt{13}}{2}$ is non-PV, this means that $\widehat{\gamma^{ }_{\mu}}$ can only have the trivial peak at zero and must essentially be continuous. 
This is in constrast with point sets arising from the Fibonacci tiling, whose diffraction measures only consist of Bragg peaks; see \cite{BGMaz}.

A more precise nature of the continuous component $\widehat{\gamma}^{ }_{\textsf{cont}}$ has only been established in \cite{BFGR} via renormalisation methods involving Lyapunov exponents, which confirmed that it is indeed purely singular continuous, that is, 
\[
\widehat{\gamma}=\dens(\varLambda)\,\delta_0+\widehat{\gamma}^{ }_{\textsf{sc}},
\] 
where $\dens(\varLambda)$ is the density of the point set $\varLambda$. 
 This is the first example of a non-PV inflation for which this is confirmed (although it has long been expected); see \cite{BGrM} for a family of non-PV inflations with the same spectral features.

For the GLB tiling mentioned in the introduction, $\lambda=\frac{5+\sqrt{5}}{2}$, which is also non-PV, and hence any derived point set cannot have non-trivial Bragg peaks (this has already been observed by Godr\`{e}che and Lan\c{c}on in \cite{GodLancon}). It was conjectured that the continuous component must be singular, which was only confirmed recently in \cite{BGM}.

We now turn to some properties of the Schr\"odinger spectrum  of this non-PV example. For any element $w\in X_{\varrho}$, one obtains the (ergodically-defined) potential 
\[
V_{w}(n)=\Phi(\sigma^n (w)), \quad\text{ for } n\in\mathbb{Z}.
\] 
Here, $\Phi\colon X_{\varrho}\to\mathbb{C}$ is usually taken as a continuous function that only depends on the letter at the origin (hyperlocal).  Using this potential function, one can define the \emph{discrete Schr\"odinger operator} on $\ell^{2}(\mathbb{Z})$ as
\[
\left[H_{w}\psi\right](n)=\psi(n+1)+\psi(n-1)+V_{w}(n)\psi(n). 
\]
One is interested in the spectrum $\Sigma_w:= \text{spec}(H_{w})$ of $H_{w}$ and the type of the spectral measures; see \cite{DF,DEG} for background and \cite{DFG} for the continuum  analogue on $L^{2}(\mathbb{R})$. 

For the non-PV example above, it follows from primitivity that  $\Sigma_w$ is constant on $X_{\varrho}$. Moreover, $\Sigma_w$ is a Cantor set of Lebesgue measure zero; see \cite{DL}. Regarding the presence of eigenvalues, note that $aaaa$ is legal since it appears in $\varrho^{3}(a)$. This, together with the hierarchical structure, allows one to use a three-block Gordon-type argument to rule out the presence of eigenvalues for almost every element in $X_{\varrho}$ (with respect to the unique invariant measure); see \cite{DEG, EG} for criteria to exclude eigenvalues.

Kotani theory and aperiodicity then rules out the presence of absolute continuous spectral measures; compare \cite{Dam,DEG}. It follows that, for a.e. $x\in X_{\varrho}$, the corresponding spectral measures are purely singular continuous.

Another well-known substitution in one dimension is the \emph{Thue--Morse} substitution
\[
\varrho^{ }_{\text{TM}}\colon \begin{cases}
a\mapsto ab\\
b\mapsto ba.
\end{cases}
\]

This is an example of a \emph{constant-length} substitution, where the symbolic and the geometric picture are equivalent. This means the underlying point sets are all supported on $\mathbb{Z}$. 
Given a weight function $w(x)\in \left\{w_a,w_b\right\}$, the diffraction $\widehat{\gamma}$ admits the form 
\[
\widehat{\gamma}=\left|\frac{w_a+w_b}{2}\right|^2\delta_{\mathbb{Z}}+\left|\frac{w_a-w_b}{2}\right|^2\widehat{\gamma}^{ }_{\text{TM}}\ast \delta_{\mathbb{Z}}. 
\]
Here, $\widehat{\gamma}^{ }_{\text{TM}}$ is called the \emph{Thue--Morse measure}, which is the first example of a purely singular continuous measure.
One gets this if one opts for the Dirac comb with \emph{balanced} weights $w_a=1, w_b=-1$ (which kills the pure point component).

We refer the reader to \cite{BGr} for a thorough discussion on the Thue--Morse measure, and \cite{BGKS, GLS, FSS} for recent developments on its multifractal features. For details regarding its $n$-point correlations, see \cite{BC}. 

On the Schr\"odinger side, the spectrum $\Sigma$ is also constant on $X_{\text{TM}}$ and is a Cantor set of Lebesgue measure zero. Kotani theory again guarantees the absence of absolutely continuous spectral measures. 

The main difference with the previous example is that ruling out eigenvalues is more subtle for Thue--Morse.
It is known that the potential associated to the \emph{palindromic} fixed point generated by the seed $a|a$ admits no eigenvalues. 
This extends to a $G_{\delta}$ subset of $X_{\text{TM}}$ via a standard orbit argument, from which one gets \emph{generic} absence of eigenvalues; see \cite{JS}.

 However, one cannot apply the same Gordon-type trick because any element $x\in X_{\text{TM}}$ is cube-free (that is, there is no legal word of the form $vvv$, where $v$ is a finite word). For this example, the \emph{almost sure} absence of eigenvalues in the Schr\"odinger spectrum is still open. 

\section{Generalisations of substitutions}

Now that we have recalled important properties of substitutions on finite alphabets and presented some examples (together with known results on their diffraction and Schr\"odinger spectra), we survey a few generalisations of substitutions. A more extensive account 
on hierarchical tilings can be found in \cite{Frank} and further examples of inflation tilings are given in \cite{Frettloeh}.  

\subsection{$S$-adic shifts}
In the iterative process above, one only uses a single rule (substitution or inflation) to generate infinite structures (either words, tilings, or point sets). One possible modification is to allow for different rules at different levels of the hierarchy. 
In the symbolic picture, this corresponds to \emph{$S$-adic directive sequences} where one gets a sequence $\boldsymbol{\varrho}:=(\varrho^{ }_0,\varrho^{ }_1,\ldots)$ of substitutions; see \cite{Berthe,Thus} for background on $S$-adic shifts.

It is worth mentioning that well-known examples of systems with pure point spectrum, namely Sturmian and Toeplitz shifts, both admit $S$-adic representations. In particular, it is known that any Sturmian shift can be generated as an $S$-adic shift using combinations of the substitutions
\[
\tau\colon \begin{cases}
0\mapsto 0\\
1\mapsto 10
\end{cases} \quad  \text{ and } \quad \phi\colon \begin{cases}
0\mapsto 01\\
1\mapsto 1;
\end{cases}
\]
see \cite{Thus}. Note that neither of the two substitutions above is primitive, but there is a notion of a primitivity for directive sequences, which allows one to get
similar properties like unique ergodicity and minimality. 
On the other hand, (one-sided) Toeplitz sequences can be obtained from directive sequences of constant-length substitutions \cite{GJ}. This freedom in mixing at different levels is not restricted to symbolic rules; see \cite{ST,SolTrevino} for extensions to geometric rules. 

\subsection{Random substitutions}
For $S$-adic shifts, the choice of the rule is determined \emph{globally} at each level (that is the rule is fixed for all letters). One can also consider the generalisation where the choice is \emph{local}, that is, for each iteration, one has the freedom where to map each symbol. Such a rule is called a \emph{set-valued substitution} or, when the choices are equipped with certain probabilities, a \emph{random substitution}; see \cite{RustSpin,Moll}. A classical example is the \emph{random Fibonacci substitution}
\[
\varrho^{ }_{\text{rand}}\colon \begin{cases}
a\mapsto &\begin{cases}
ab,& \text{with probability }p\\
ba,& \text{with probability }1-p
\end{cases}\\
b\mapsto & a;
\end{cases}
\]
see \cite{RustSpin,BSS} for its properties and \cite{BaakeMoll,EMM,Moll} for generalisations (including the random noble means and their counterparts in larger alphabets).

Here, one has a choice to map $a$ to either $ab$ or $ba$ at each step. Subshifts associated to random substitutions are highly non-minimal, admit many ergodic measures \cite{GS}, and have positive entropy \cite{GMRS,EMM}.
These subshifts realise a wide range of growth behaviour for the complexity function \cite{Mitchell}.  
They can also contain (shift)-periodic points \cite{Rust}. 
As in the case of $S$-adic shifts, there is also a notion of primitivity for random substitutions. 

Under certain geometric compatibility conditions, one can also derive tilings and point sets from random substutitions. Such systems can be considered hybrids of random tilings (where there is absolutely no hierarchical structure present) and substitutions (where everything is deterministic). 

 The random Fibonacci example satisfies such a geometric requirement, and one can unambiguously assign a tile of length $\tau=\frac{1+\sqrt{5}}{2}$ to $a$ and a tile of length $1$ to $b$.

Whenever $p>0$, the (almost sure) diffraction $\widehat{\gamma}$ in this case has a non-trivial absolutely continuous component due to the presence of randomness. More precisely, 
\[
\widehat{\gamma}=\sum_{k\in \frac{\mathbb{Z}[\tau]}{\sqrt{5}}}I_p(k)\delta_k+\phi_p(x)\mu^{ }_{\text{Leb}}, 
\]
where closed forms for both the (probability-dependent) intensities $I_p(k)$ for the pure point part and the Radon--Nikodym density $\phi_p(x)$ for the absolutely continuous component are available; see \cite{BSS} for details.

\subsection{Multiscale substitutions}

The self-similar nature of the hierarchical structure in classical substitutions manifests in the presence of a single inflation factor (which is also called a scaling constant). Supertiles at level $n$ are (uniformly related) to supertiles at level $n{+}1$ by this scale $\lambda$ due to the inflate-and-subdivide rule. 

\begin{figure}[!h]
\[
\underset{\mathfrak{t}}{\includegraphics[scale=1.0,valign=c]{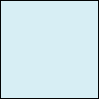}}
 \,\,\,\overset{\lambda}{\longmapsto} \,\,\,
\begin{cases}
\hspace{3mm}
\includegraphics[scale=1.0,valign=c]{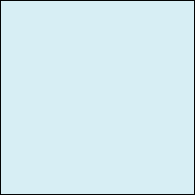}\,, & \text{if }\lambda\cdot \text{vol}(\mathfrak{t})<1, \vspace{5mm}\\
\hspace{3mm} \includegraphics[scale=1.0,valign=c]{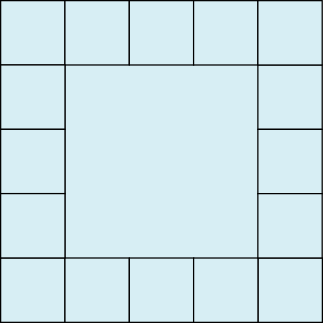}\,,& \text{otherwise}. 
\end{cases}
\]
\caption{An example of a multiscale substitution with a square prototile.}\label{fig:multiscale}
\end{figure}

A generalisation of this notion is the introduction of multiple scales by tweaking the 
``inflate-and-subdivide" concept into an ``inflate-and-only subdivide when a certain threshold is reached". Tilings generated by such a rule are called \emph{multiscale subtitution tilings}. An example in Figure~\ref{fig:multiscale} is taken from \cite{SS}, where the image of a square under inflation is only subdivided once the area of the blown-up square exceeds or equals $1$. 

Here the scales are the distinct (length scales) which appeear in the subdivided 
image of a tile with unit volume. For the example in Figure~\ref{fig:multiscale}, 
there are \emph{two scales}, namely $\frac{1}{5}$ and $\frac{3}{5}$, since there are squares with these side lengths appearing in a subdivided unit square.

In this setting, one advantage is that one goes beyond the discrete steps (levels) of substitution and one can consider the evolution of a tile under a continuous process of being substituted. This yields a continuous collection of supertiles controlled by a time parameter. 
We refer to \cite{SS,Smilansky,Smilansky-2} for detailed accounts on multiscale substitutions, \cite{Frank-Sadun2,Sadun-2} for predecessors, and \cite{AF} for an account on the related \emph{interval splitting procedure} due to Kakutani.

\subsection{Fusion tilings}
Another family of hierarchical rules is that of \emph{fusion rules}, where the hierarchies are not constructed by the usual inflating procedure. Rather, one constructs a level-$n$ patch $P_n(i)$ by putting together level-$(n{-}1)$ patches $P_{n-1}(j)$, where $i,j$ are prototile labels; see \cite{Frank-Sadun}. 
Note that this in general does not give a well-defined inflation factor, but the substitution matrix finds a natural analogue in the \emph{transition map} $M_{n,n-1}$, which contains information on the statistics of patches of certain type. In the case of an FLC fusion rule, one has
\[
\left(M_{n,n-1}\right)^{ }_{ij}= \text{number of } P_{n-1}(i) \text{ in } P_n(j).
\]
Note that all FLC tilings are fusion tilings. There is also a formalism for tilings with infinite local complexity \cite{Frank-Sadun2}.

This flexibility allows one to construct lots of interesting examples (in view of dynamical, geometric, and spectral properties). In particular, the \emph{scrambled Fibonacci} tilings, which are generated by modifying the level-$n$ supertiles of the Fibonacci inflation rule by  introducing defects at sparse levels yield repetitive, FLC point sets which are \emph{not Meyer} but still have pure point diffraction. Moreover, none of the Bragg peaks are topological (that is, no dynamical eigenvalue admits a continuous eigenfunction); see \cite{Frank-Sadun2,KS} for details.

\subsection{Non-Abelian generalisations}

All constructions we have seen so far give rise to tilings of $\mathbb{Z}^d$ or $\mathbb{R}^d$, which are Abelian groups. Here we would like to mention recent extensions to other spaces. In \cite{BHP}, the idea of symbolic substitutions is extended to homogeneous spaces, such as Lie groups. One can extend the notion of repetitivity, this time with respect to the natural metric of the underlying space; see \cite{BHP} for an example of a substitution on the Heisenberg group $\mathbb{H}_3(\mathbb{R})$. One can also consider substitutions on (coloured) infinite trees; see \cite{BL} for a generalisation of the Thue--Morse substitution on the (free) semigroup $\mathbb{F}^{+}_2$ on two generators.

\section{Substitutions on compact alphabets}\label{sec:compact}

Let us now move to the main focus of this review, which are substitutions on compact alphabets. We first mention the differences compared to the generalisations mentioned in the previous section. In what follows, we only consider a rule generated by a \emph{single  substitution} (and not a sequence thereof as in $S$-adic shifts). All the rules we present are \emph{deterministic} and no randomness is involved. Whenever present, the inflation rules they generate only have a \emph{single scale}. Lastly, they can be considered as fusion rules, where the transition map is constant for all $n$.

We start with the symbolic treatment as in Section~\ref{sec:1d-finite}. Here, our alphabet $\mathcal{A}$ is a \emph{compact} Hausdorff topological space. This means the set $\mathcal{A}^{+}$ of finite words and the set $\mathcal{A}^{\mathbb{Z}}$ of all bi-infinite words over $\mathcal{A}$ are both compact spaces as well. 
Now, a substitution $\varrho$ on $\mathcal{A}$ is a \emph{continuous} rule which maps letters in $\mathcal{A}$ into finite words over $\mathcal{A}$.

This continuity condition on $\varrho$ is trivially satisfied in the case of finite alphabets. In the compact alphabet setting, we will see that it has lots of strong consequences. For example, if $\mathcal{A}$ is connected, this immediately implies that any substitution on $\mathcal{A}$ must be of constant length.

 The series of works \cite{MRW,MRW-2,FGM,FGM-2} of the author and his collaborators investigate geometric, combinatorial, dynamical, operator-theoretic, and spectral properties of these objects.

Before we proceed, we note that substitutions on infinite alphabets have been investigated before; see \cite{RY, AS-2, Queffelec} for their appearances in the context of automatic sequences and their generalisations. Those on countable alphabets also remain an active area of research, with their links to countable Markov chains; see \cite{Ferenczi,BJS}. Finally, dynamical systems coming from substitution rules on compact spaces have already appeared in \cite{DOP,Frank-Sadun2}.

From $\varrho$, we want to build a subshift $X_{\varrho}\subset \mathcal{A}^{\mathbb{Z}}$. In general, $\varrho$ does not admit a fixed point. Nevertheless, we can use the \emph{language} $\mathcal{L}_{\varrho}$ of $\varrho$ to build $X_{\varrho}$. The language is defined as
\[
\mathcal{L}_{\varrho}=\overline{\{v\in \mathcal{A}^{+}\colon v \text{ appears in some } \varrho^{n}(a)\} },
\]
where the closure is taken with respect to the topology of $\mathcal{A}^{+}$.
The subshift is now defined as 
\[
X_{\varrho}=\left\{x\in \mathcal{A}^{\mathbb{Z}}\colon \text{if }v\text{ is a subword of } x \implies v\in \mathcal{L}_{\varrho}\right\}.
\]
As before, if we equip it with the left shift action, $(X_{\varrho},\sigma)$ becomes a topological dynamical system.

Let us give two guiding examples which we will revisit as we go on. 
Let $\mathcal{A}=\mathbb{N}_0\cup\left\{\infty\right\}$ be the one-point compactification of the nonnegative integers. Consider the rule 
\begin{equation}\label{eq:subs-Nstar}
\varrho^{ }_{\infty}\colon \begin{cases}
0\mapsto& 0\,1\\
1\mapsto& 0\, 0\, 1\\
\vdots&\\
n\mapsto & 0\, (n{-}1)\,(n{+}1)
\end{cases}
\end{equation}
Since we want $\varrho$ to be continuous, we immediately get $\infty\mapsto 0\,\infty\,\infty$.

Another example is a generalisation of the Thue--Morse substitution on $\mathcal{A}=\mathbb{T}\simeq \mathbb{R}/\mathbb{Z}$ given by
\begin{equation}\label{eq:subs-S1}
\varrho^{ }_{\alpha}\colon (\theta)\mapsto (\theta)\, (\theta{+}\alpha)\quad\text{for }\theta\in \mathbb{T},
\end{equation}
where $\alpha\in \mathbb{T}$ is fixed.

In the compact alphabet case, \emph{primitivity} means that, given a non-empty open set $U\subset\mathcal{A}$, one can find a power $n:=n(U)$ such that 
\[
\varrho^{n}(a)\text{ contains a letter in }U, \text{ for all }a\in\mathcal{A}. 
\]
Going back to the examples above, one can see that $\varrho^{ }_{\infty}$ is primitive because every level-$1$ superword starts with a $0$ and $\varrho_{\infty}^n(0)$ has the letter $n$ as a subword.
Thus, if $U=\left\{m\right\}$, one can pick $n(U)=m+1$. Otherwise, if $U$ is a ball around $\infty$, one only needs to take a power $n$ such that $n{-}1$ is in $U$.  
For the generalised Thue--Morse, it is easy to see via a density argument that $\varrho^{ }_{\alpha}$ is primitive if and only if $\alpha$ is irrational. 

As in the finite alphabet case, primitivity implies that $(X_{\varrho},\sigma)$ is minimal and that every element $x\in X_{\varrho}$ is \emph{almost repetitive} (where finite words appear with bounded gaps if one allows some local perturbations); compare \cite{FR} for the notion of almost repetitivity for Delone sets. 

We now discuss some operator-theoretic aspects of $\varrho$. 
To $\mathcal{A}$, one associates the space $C(\mathcal{A})$ of (real-valued) continuous functions on $\mathcal{A}$. This is a Banach space (with the norm $\|f\|=\max_{a\in\mathcal{A}}|f(a)|$) and is indeed a Riesz space ($V$-lattice), making it a \emph{Banach lattice}. Crucial to our analysis is the 
positive cone
\[
K=\left\{f\in C(\mathcal{A})\colon f(a)\geqslant 0\text{ for all }a\in \mathcal{A}\right\}.
\] 

The \emph{substitution operator} $M\colon C(\mathcal{A})\to C(\mathcal{A})$ (which is the generalisation of the transpose of the substitution matrix) is defined via
\[
(Mf)(a)=\sum_{b\in\varrho(a)} f(b), 
\]
where the sum is taken over all letters $b$ (counted with multiplicity) that appear in $\varrho(a)$. As an example, if we consider the substitution operator of $\varrho^{ }_\infty$, we get $(Mf)(\infty)=f(0)+2f(\infty)$, for any $f\in C(\mathcal{A})$. Since $MK\subset K$, $M$ is a \emph{positive} operator. We refer the reader to the monographs  \cite{EFHN, SchaeferBook} and to \cite{Karlin} for properties of positive operators on Banach lattices.

\begin{figure}[!h]
\begin{center}
\includegraphics[scale=0.8]{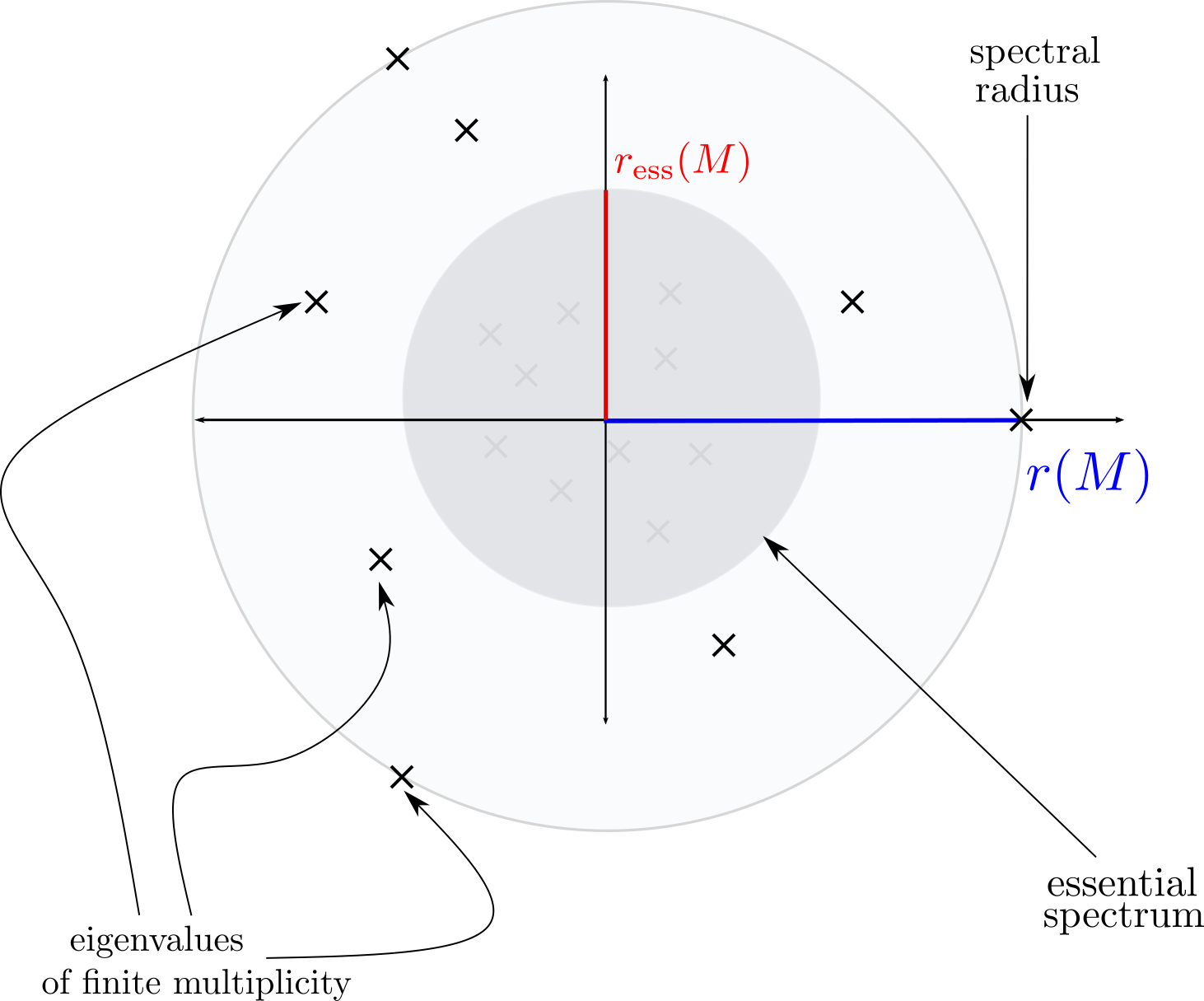}
\end{center}
\caption{Illustration of the spectrum $\text{spec}(M)$ of a (positive) quasicompact operator $M$.}\label{fig:QC}
\end{figure}

Due to positivity,  the \emph{spectral radius} $r(M)$ of $M$ is in the spectrum 
\[\text{spec}(M)=\left\{\lambda\in \mathbb{C}\colon (\lambda\mathbb{I}-M) \text{ is not invertible}\right\}. \]
We call $M$ \emph{quasi-compact} if there is a power $n\in\mathbb{N}$ and a compact operator $C$ such that 
\[
\|M^n-C\|_{\text{op}}<r(M)^n\quad. 
\]
Equivalently, $r_{\text{ess}}(M)<r(M)$
where $r_{\text{ess}}(M)$ is the essential spectral radius of $M$. Outside the essential spectrum, there are only finitely many elements of $\text{spec}(M)$, which are all eigenvalues of finite multiplicity; see Figure~\ref{fig:QC}. For a succinct account on quasi-compact operators, one can consult \cite{HH}.

We call $\boldsymbol{L}\in K\subset C(\mathcal{A})$ a \emph{natural length function} if there exists $\lambda >0$ such that $M\boldsymbol{L}=\lambda \boldsymbol{L}$ (that is $\boldsymbol{L}$ is a right eigenvector of $M$). 
Note that the existence of such a function with $\lambda >1$ allows one to derive an inflation rule from $\varrho$, as in Section~\ref{sec:1d-finite}.

In the compact alphabet setting, we then have the following questions regarding a substitution $\varrho$

\begin{enumerate}
\item When does $\varrho$ admit a natural length function?
\item When is $(X,\varrho)$ uniquely ergodic?
\end{enumerate}

Note that primitivity alone does not suffice for (2); see \cite{DOP} for a counterexample. This means that we need an additional assumption, which is given in the following result.

\begin{theorem}[\cite{MRW}]\label{thm:main-result}
Let $\varrho$ be substitution on a compact Hausdorff alphabet $\mathcal{A}$. If $\varrho$ is primitive and $M$ is quasicompact, then, 
\begin{enumerate}
\item $\varrho$ admits a strictly positive natural length function $\boldsymbol{L}$ with $\lambda=r(M)$ and
\item $(X_\varrho,\sigma)$ is uniquely ergodic.  \qed
\end{enumerate} 
\end{theorem}

Note that $M$ is always compact (and hence quasicompact) in the finite alphabet setting, and hence the extra assumption in the result above is trivially satisfied. 
The sufficient conditions presented above are not the weakest assumptions, but they are normally easier to confirm for examples; see \cite{MRW} for a combinatorial criterion for quasi-compactness. 

We briefly comment on the implications. For (1), we even have something more: this length function is unique in a strong sense: it is the unique eigenvector to the eigenvalue $r(M)$ and it is the unique eigenvector of $M$ which is strictly positive (both up to scaling). 
 For (2), we also get that the induced tiling dynamical system $(\Omega,\mathbb{R})$ is uniquely ergodic. 

The tilings in $\Omega$ whose tiles are defined by $\boldsymbol{L}$ 
are typically ILC, but still define Delone point sets because the tile lengths satisfy
\[
0<\ell^{ }_{\min}\leqslant\boldsymbol{L}(a)\leqslant \ell^{ }_{\max},
\]
for some constants $\ell^{ }_{\min},\ell^{ }_{\max}$, by continuity. In particular, one can get Delone sets with inflation symmetry but which are not of finite type, generalising those in \cite{Lagarias2}.

\section{Examples}\label{sec:examples}

\subsection{Substitutions on compactifications of $\mathbb{N}$}

The substitution $\varrho^{ }_{\infty}$ in Eq.~\eqref{eq:subs-Nstar} satisfies the conditions of Theorem~\ref{thm:main-result}. One can confirm the quasi-compactness of $M$ via a checkable criterion. 

For $\varrho^{ }_{\infty}$, the inflation factor is $\lambda=r(M)=5/2$. This is not possible for finite alphabet substitutions, where $\lambda$ is always an algebraic integer, which $5/2$ is not. The corresponding tile lengths (with normalisation $\boldsymbol{L}(0)=1$) are given by 
\[
\boldsymbol{L}(n)=2-\frac{1}{2^n},\quad \quad  \boldsymbol{L}(\infty)=2,
\]
while the (uniform) relative letter frequencies are 
\[
\nu(n)=\frac{1}{2^{n+1}},\quad\quad \nu(\infty)=0,
\]
for any $x\in X_{\varrho}$; see \cite{MRW-2}.

The example in Eq.~\eqref{eq:subs-Nstar} can be generalised using a sequence 
$\boldsymbol{m}=(m_i)_{i\geqslant 0}$ of non-negative integers as follows. For $n\in\mathbb{N}_{0}$, one builds the rule
\[
\varrho^{ }_{\boldsymbol{m}}\colon \begin{cases}
0\mapsto& 0^{m_0}\,1\\
1\mapsto& 0^{m_1}\, 0\, 1\\
\vdots&\\
n\mapsto & 0^{m_n}\, (n{-}1)\,(n{+}1),
\end{cases}
\]
where, $0^{m}$ refers to the concatenation of $m$ copies of $0$. In particular, the constant sequence $\boldsymbol{m}=(1)_{i\geqslant 0}$ corresponds to $\varrho^{ }_{\infty}$.

One then considers an appropriate embedding $\iota(\mathbb{N}_0)$ of $\mathbb{N}_0$ in some full shift, where the closure is taken to give the alphabet $\mathcal{A}=\overline{\iota(\mathbb{N}_0)}$.
Both choices of the ambient full shift and the embedding depend on the sequence $\boldsymbol{m}$. 
Also depending on $\boldsymbol{m}$, the set $\mathcal{A}\setminus \iota(\mathbb{N}_0)$ of accumulation points can be finite (as in $\varrho^{ }_{\infty}$ in Eq.~\eqref{eq:subs-Nstar}), countably infinite, or even uncountable; see \cite{FGM} for details. 

\begin{figure}[!h]
\begin{center}
\includegraphics[scale=3.0]{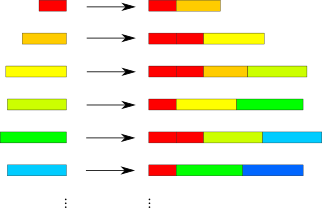}
\caption{The first few supertiles of the inflation rule $\varrho^{ }_{\boldsymbol{m}} $ with transcendental $\lambda$ corresponding to $\boldsymbol{m}=(1,2,2,1,2,1,1,2,2,1,1,2,\ldots)$. The  lengths are normalised so that the tile $\mathfrak{t}_0$ has length $1$ (red). Taken from the Tilings Encyclopedia \cite{FGH}.}\label{fig:trans}
\end{center}
\end{figure}

This construction is used to realise \emph{any} $\lambda>2$ as an inflation factor of a substitution on a compact alphabet.
Under some natural assumptions on $\boldsymbol{m}$, it was shown that $\varrho^{ }_{\boldsymbol{m}}$ satisfies the conditions of Theorem~\ref{thm:main-result}; see \cite{FGM} for a complete account.

Apart from this, $\lambda$ and $\boldsymbol{L}(n)$ (for $n\in\mathbb{N}_0$) admit the closed forms  
\[
\lambda=\mu+\frac{1}{\mu},\quad \quad \quad\boldsymbol{L}(n)=\mu^n+\sum_{j=1}^n\sum_{i=j}^{\infty} m_i\mu^{i+n+1-2j},
\]
where $\mu$ is the unique solution in the interval $(0,1)$ of the equation, $\frac{1}{\mu}=\sum_{i=0}^{\infty}m_i\mu^i$. This allows one to estimate the inflation factor and the lengths numerically to the desired precision. The relative letter frequencies 
(for $n\in \mathbb{N}_0$) are given by $\nu(n)=(1-\mu)\mu^{n}$. Any accumulation point has relative frequency zero.

One interesting member of this (infinite) family of substitutions is the case when $\boldsymbol{m}$ is the Thue--Morse sequence on $\{1,2\}$, which is given by 
\[
\boldsymbol{m}=(1,2,2,1,2,1,1,2,2,1,1,2,\ldots)
\]
Because $\boldsymbol{m}$ is non-periodic and repetitive, $\mathcal{A}\setminus \iota(\mathbb{N}_0)$ is uncountable. 

Using the same techniques as in \cite{FGM}, one can show that the corresponding $\lambda$ for this explicit choice of $\boldsymbol{m}$ is a \emph{transcendental} number. This means $\lambda$ is not a root of any (non-zero) polynomial whose coefficients are in $\mathbb{Q}$. 
Here, $\lambda\approx 2.6113$. The images under the inflation rule for the first few tiles $\mathfrak{t}_0$ up to $\mathfrak{t}_5$ are given in Figure~\ref{fig:trans}.

\subsection{Substitutions on $\mathbb{T}$}
We now turn to the family of substitutions $\varrho^{ }_{\alpha}$ in Eq.~\eqref{eq:subs-S1}, with $\alpha\in [0,1)$ irrational. They are all constant-length so they trivially have $\lambda=2$ and the length function $\boldsymbol{L}(a)=1$ for all $a\in\mathbb{T}$, so the symbolic and geometric pictures agree.

In this case, the substitution operator $M$ is never quasicompact. However, $M$ is \emph{strongly power convergent}, that is $\left(M/r(M)\right)^n(f)\rightarrow Pf$
for some bounded operator $P$, for all $f\in C(\mathcal{A})$. It turns out that this, together with the primitivity of $\varrho^{ }_{\alpha}$, is enough to get unique ergodicity; see \cite{MRW}.  

As in the case of the Thue--Morse substitution, all coloured point sets, and hence all weighted Dirac combs, are supported in $\mathbb{Z}$, which is an advantage in computing the autocorrelation $\gamma$. Here, it admits the form \[
\gamma=\sum_{z\in \mathbb{Z}}\eta(z)\delta_z
\] 
where the \emph{autocorrelation coefficients}
\[
\eta(z)=\lim_{n\to \infty} \frac{1}{2n+1}\sum_{x,x-z\in [-n,n]}w(x)\overline{w(x-z)}
\]
exist uniformly on $X_{\varrho}$ because of unique ergodicity. In what follows, we fix the weight function $w(x)=\ee^{2\pi \ii\theta}$ if the point at position $x$ is of type $\theta$. Note that this choice of weight is balanced since $\int_{\mathbb{T}} \ee^{2\pi \ii\theta}\dd\mu^{ }_{\text{H}}(\theta)=0$, where $\mu^{ }_{\text{H}}$ is the Haar measure on $\mathbb{T}$.   

\begin{figure}[!h]
\includegraphics[scale=0.85]{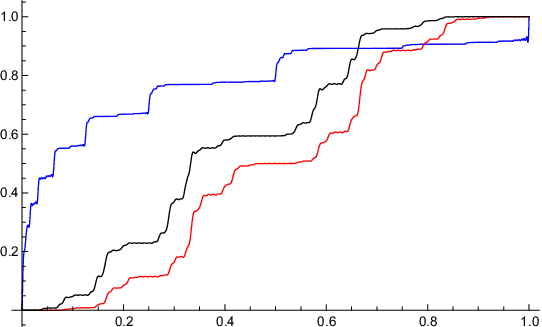}
\caption{Plots of the approximation of the distribution function $F_{\alpha}(x):=\widehat{\gamma}^{ }_\alpha([0,x))$ for $\alpha=1/2$ (red), $\alpha=\left\{\sqrt{2}\right\}$ (black), and $\alpha=\left\{\pi\right\}$ (blue). The choice $\alpha=1/2$ corresponds to the classical Thue--Morse measure.}
\label{fig:distrib}
\end{figure}

For a fixed $\alpha$, $\eta(z)$ satisfies the recursive relations 
\begin{align*}
\eta(2z)&=\eta(z)\\
\eta(2z+1)&=\frac{1}{2}\left(\eta(z)\ee^{-2\pi \ii\alpha}+\eta(z+1)\ee^{2\pi \ii\alpha}\right)
\end{align*}
for $z\in \mathbb{Z}$. Since $\eta(0)=1$, one gets 
\[
\eta(1)=\frac{\ee^{-2\pi\ii\alpha}}{2-\ee^{2\pi\ii\alpha}}.
\]
In \cite{MRW-2}, it was shown that, for a family of constant-length substitutions containing $\varrho^{ }_{\alpha}$, $|\eta(1)|<1$ implies that the diffraction $\widehat{\gamma}$ is purely singular continuous. 
Since the weighted point set is supported on $\mathbb{Z}$, one can write the diffraction as $\widehat{\gamma}=\delta_{\mathbb{Z}}\ast \widehat{\gamma}^{ }_{\alpha}$. 
One can easily check that $|\eta(1)|=1$ for the example above only when $\alpha=0$. This means that, for all non-zero $\alpha$, the corresponding diffraction is purely singular continuous. Moreover, the $\mathbb{Z}$-periodic part $\widehat{\gamma}^{ }_{\alpha}$ has a representation as a Riesz product given by 
\[
\widehat{\gamma}^{ }_{\alpha}=\prod_{n=0}^{\infty}\frac{1}{2}\left|1+\ee^{2\pi\ii(2^{n}k+\alpha)}\right|^2
\]
These Riesz products arising from generalised Thue--Morse polynomials are studied in the context of multifractal analysis in \cite{FSS}. 
The plots of the  distribution functions for $\alpha=\left\{1/2,\{\sqrt{2}\},\{\pi\}\right\}$ are provided in Figure~\ref{fig:distrib}.  Here, $\{y\}$ denotes the fractional part of $y$. 
  
The substitution $\varrho^{ }_{\alpha}$ can be further generalised into a block substitution in two dimensions. Here, a prototile is a square tile of side length $1$ labelled by a letter $\theta\in \mathbb{T}$. An example is given by the rule\\
\vspace{2mm}

$\varrho\,\,\colon$
\begin{squarecells}{1}
$\theta$  \nl
\end{squarecells}\quad$\mapsto$\quad 
\begin{squarecells}{3}
$\theta$ & $\theta{+}\frac{\sqrt{7}}{2}$ & $\theta$ \nl
$\theta$ & $\theta{+}\frac{\sqrt{3}}{2}$ & $\theta{+}\frac{1}{2}$ \nl
$\theta$ & $\theta{+}\frac{1}{4}$ & $\theta{+}\frac{3}{4}$ \nl
\end{squarecells}\\
\vspace{2mm}

\begin{figure*}[!h]
\captionsetup{width=.9\linewidth}
\begin{center}
\includegraphics[scale=0.6]{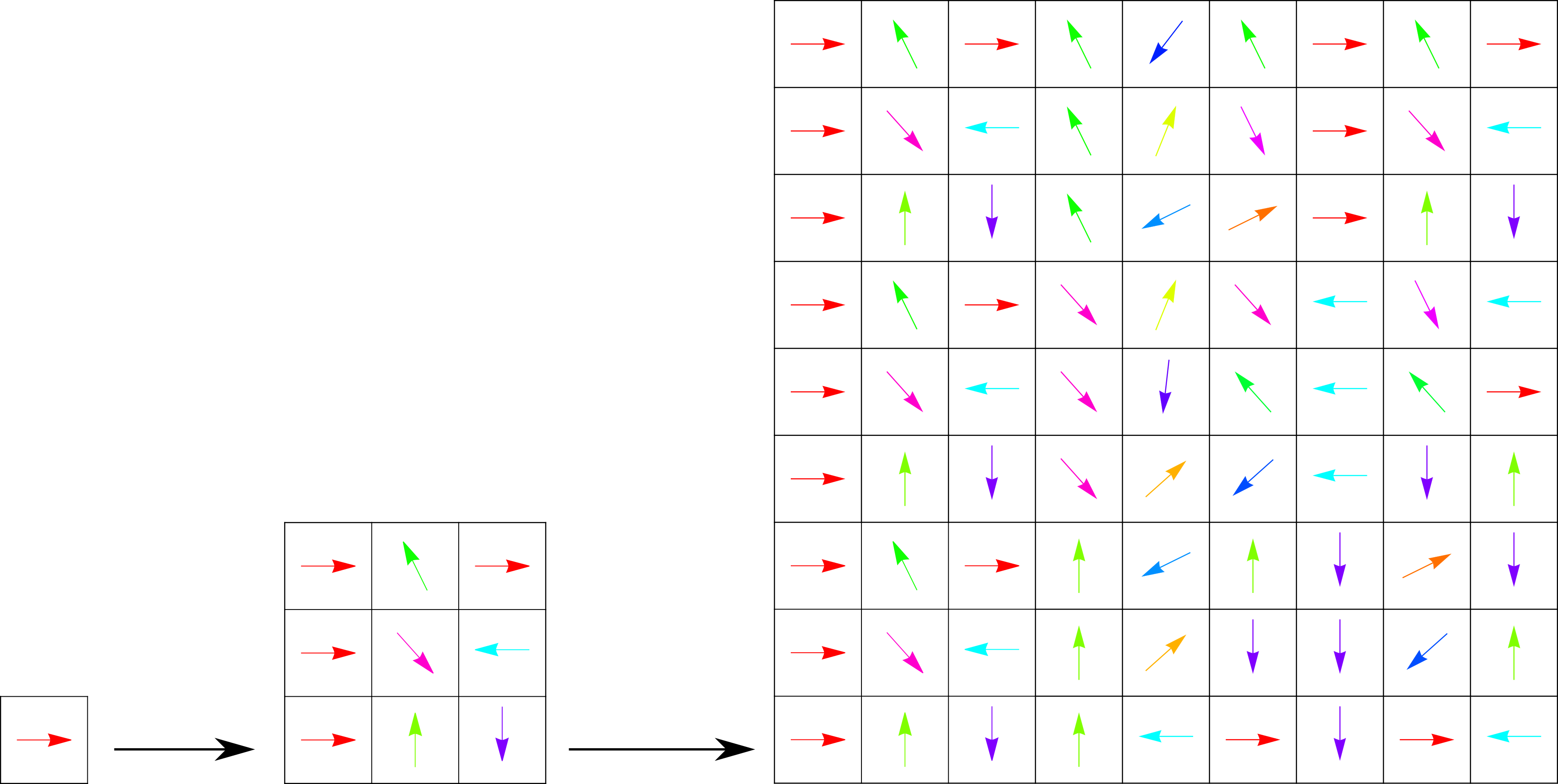}
\end{center}
\caption{The first three supertiles for the planar substitution on $\mathcal{A}=\mathbb{T}\simeq [0,1)$ for $\theta=0$. Here, we represent the label by an arrow pointing in the direction of the weight $w(\theta)=\ee^{2\pi\ii\theta}\in S^1$ we used in computing the diffraction. To aid visualisation, the arrows are also coloured following a gradient with $0$=red, $\frac{1}{4}$=green, $\frac{1}{2}$=blue, and $\frac{3}{4}$=purple.}\label{fig:bij}
\end{figure*}

The level-$0$,$1$, and $2$ supertiles corresponding to $\theta=0$ are shown in Figure~\ref{fig:bij}.
The presence of irrational rotations guarantees the primitivity of  $\varrho$. Moreover, the two-dimensional system $(X_{\varrho},\mathbb{Z}^2)$ is uniquely ergodic.  Here, to get the singular continuity of $\widehat{\gamma}$, one needs to check whether $|\eta(e_i)|<1$, where $\left\{e_i\right\}$ is standard basis for $\mathbb{Z}^2$.  Using the same techniques as in one dimension, one can show that the weighted point set (with the same weight function as above) from any element of $X_{\varrho}$ has purely singular continuous diffraction; compare \cite{MRW-2}.

\section{Summary and Outlook}\label{sec:outlook}
Theorem~\ref{thm:main-result} provides conditions for geometric realisability and for uniformity of statistical properties of substitutions on compact alphabets. 
As we saw in the examples presented, these already revealed phenomena which are not 
possible in the finite alphabet case. The next step is to look closely at further spectral implications of these results. As an example, one can consider Schr\"odinger potentials generated by these subshifts. These potentials are ergodically-defined, can take infinitely many values (but still in a compact subset of $\mathbb{C}$), are generally not (purely) quasi-periodic, yet still admit a hierarhical structure.  As an outlook, we want to uncover further consequences of having this structure, which we plan to address in a series of future works.

% to add a large picture spanning two columns: 
%\begin{figure*}
%
%\end{figure*}

	\section*{Acknowledgements}
	The author would like to explicitly mention that the methods, results, and examples presented in Sections~\ref{sec:compact} and ~\ref{sec:examples} of this review are based on a series of joint works with Dan Rust (Milton Keynes, UK) and Jamie Walton (Nottingham, UK), and with Dirk Frettl\"oh (Bielefeld, Germany) and Alexey Garber (Brownsville, TX, USA). The author would like to thank Michael Baake, Michael Coons, Jake Fillman, Natalie Frank, Philipp Gohlke, Lorenzo Sadun, and Nicolae Strungaru for discussions, and Jan Maz\'{a}\v{c} and Andrew Mitchell for comments on the manuscript.  This work was supported by the German Research Foundation (Deutsche Forschungsgemeinschaft, DFG) through SFB 1283/2 2021–317210226.
	\bigskip

	\bibliographystyle{amsart}

\end{document}